\documentclass{article}
\usepackage{amsmath,amssymb,latexsym,a4}
\usepackage[title]{appendix}

\newtheorem{theorem}{Theorem}
\newtheorem{lemma}{Lemma}

\newtheorem{corollary}[lemma]{Corollary}
\newtheorem{proposition}[lemma]{Proposition}

\newtheorem{definition}{Definition}
\newtheorem{example}[lemma]{Example}
\newtheorem{remark}[lemma]{Remark}

\newcommand{\bl}{\begin{lemma}}
\newcommand{\el}{\end{lemma}}
\newcommand{\bt}{\begin{theorem}}
\newcommand{\et}{\end{theorem}}
\newcommand{\bcor}{\begin{corollary}}
\newcommand{\ecor}{\end{corollary}}
\newcommand{\bp}{\proof{.}}
\newcommand{\ep}{\eop}
\newcommand{\bpr}{\begin{proposition}}
\newcommand{\epr}{\end{proposition}}
\newcommand{\brem}{\begin{remark} \em}
\newcommand{\erem}{\end{remark}}
\newcommand{\bd}{\begin{definition} \em}
\newcommand{\ed}{\end{definition}}
\newcommand{\bex}{\begin{example} \em
}
\newcommand{\eex}{\end{example}}
\newcommand{\beq}{\begin{equation} }
\newcommand{\eeq}{\end{equation}}

\newcommand{\bi}{\begin{itemize}
  }
\newcommand{\ei}{\end{itemize}}
\newcommand{\ben}{\begin{enumerate} }
\newcommand{\een}{\end{enumerate} }

\newcommand{\refeq}[1]{(\ref{#1})}
\newenvironment{enumr}{

\begin{enumerate}     }{\end{enumerate}

}

\newcommand{\benr}{\begin{enumr}
  }
\newcommand{\eenr}{
\end{enumr}}

\newcommand{\ignore}[1]{}

\newcommand{\al}[1]{\forall #1\:}

\newlength{\hilflh}

\newcommand{\cL}{{\mathcal L}}

\newcommand{\cF}{{\mathcal F}}
\newcommand{\cE}{{\mathcal E}}

\newcommand{\cP}{{\mathcal P}}

\newcommand{\fB}{\mathfrak{B}}

\newcommand{\fI}{\mathfrak{I}}

\newcommand{\ga}{\alpha}
\newcommand{\gb}{\beta}
\newcommand{\gd}{\delta}
\renewcommand{\ge}{\varepsilon}

\newcommand{\gs}{\sigma}
\newcommand{\gy}{\gamma}
\newcommand{\gw}{\omega}
\newcommand{\gS}{\Sigma}

\renewcommand{\phi}{\varphi}

\newcommand{\rst}{\upharpoonright}

\newcommand{\fc}{\Vdash}      

\newcommand{\Lim}{\mathrm{Lim}}

\newcommand{\iffdef}{\stackrel{\text{def}}{\iff}}

\renewcommand{\leq}{\leqslant}
\renewcommand{\geq}{\geqslant}

\newcommand{\eop}{$\Box$ \protect\par \addvspace{\topsep}}
\newcommand{\proof}[1]{\protect\par\addvspace{\topsep}\noindent {\bf Proof#1}}

\newcommand{\Rcn}{{\mathrm{RC}^\nab}}

\newcommand{\nab}{\nabla\hspace{-0.7pt}}

\begin{document}

\title{A universal Kripke frame for the variable-free fragment of $\Rcn$}

\author{Lev D. Beklemishev\thanks{Research financed by a grant of the Russian Science Foundation (project No. 16-11-10252).} \\
Steklov Mathematical Institute,
RAS \\ email: \texttt{bekl@mi.ras.ru}}

\maketitle

\begin{abstract}
This note extends \cite{Bek18a} by characterizing a universal Kripke frame for the variable-free fragment of the reflection calculus with conservativity operators $\Rcn$. The frame here is obtained from the set of all filters on the Ignatiev $\Rcn$-algebra which is an isomorphic presentation of the Lindenbaum--Tarski algebra of the variable-free fragment of $\Rcn$. We give a constructive `coordinatewise' characterization of the set of filters and of the frame relations corresponding to the modalities of the algebra.
\end{abstract}

In view of Theorem 3 of \cite{Bek18a} it is natural to ask if one can describe a universal Kripke frame for the variable-free fragment of $\Rcn$. There is a general construction associating with a SLO $\fB=(B,\land^\fB,\{a^\fB:a\in\Sigma\})$ its `dual' Kripke frame, similar to the way the canonical model of an strictly positive logic $L$ is obtained from its Lindenbaum--Tarski algebra.

Recall that a \emph{filter in $\fB$} is a nonempty subset $F\subseteq B$ such that
\ben\item If $x\leq_\fB y$ and $x\in F$ then $y\in F$;
\item If $x,y\in F$ then $x\land^\fB y\in F$.
\een
The set of filters of $\fB$ will be denoted $\cF(\fB)$. On $\cF(\fB)$ one can define binary relations $\{R_a:a\in \Sigma\}$ as follows: For all $F,G\in \cF(\fB)$,
$$
F R_a G \iffdef \al{x\in G} a^\fB(x)\in F.
$$
Let $\fB^*$ denote the Kripke frame $(\cF(\fB),\{R_a:a\in\Sigma\})$ together with the canonical valuation $v:\fB\to \cP(\cF(\fB))$, where $v(x):=\{F\in \cF(\fB): x\in F\}$.
It is then easy to see that, for all $x\in \fB$ and $a\in\gS$, $R_a^{-1}(v(x))=v(a^{\fB}(x))$. Hence, we obtain the following corollaries.

\bpr \benr \item The map $v:\fB\to \cP(\fB^*)$ is an embedding of $\fB$ into the algebra $(\cP(\fB^*),\cap,\{R^{-1}_a:a\in\Sigma\})$.
\item If $A,B$ in $\cL_\gS$ are variable-free, then $A\vdash B$ holds in $\fB$ iff $\fB^*,F\fc A\to B$ for all $F\in \fB^*$.
    \eenr
\epr

\bcor $\fI^*$ is complete for the variable-free fragment of $\Rcn$. \ecor

Now our task is to give a more explicit description of the Kripke frame $\fI^*$.
To this end we analyse the structure of filters of $\fI$.
Given a subset $F\subseteq \fI$ let $F_i:=\{\ga_i:\vec\ga\in F\}$ denote the $i$-th projection of $F$.
\bl \label{f-clos} If $F$ is a filter in $\fI$, then the following conditions hold:
\benr
\item Each $F_i$ is nonempty;
\item If $\ga\in F_i$ and $\gb<\ga$ then $\gb\in F_i$;
\item If $\ga\in F_i$, $\gb\in F_{i+1}$ and $\ell(\ga)<\gb$ then $\ga+\gw^\gb\in F_i$.
\eenr
\el

\bp\ Claims (i) and (ii) are easy, we prove claim (iii). Let $\vec\ga,\vec\gb\in F$ such that $\ga_i=\ga$ and $\gb_{i+1}=\gb$. Then $\vec\gy=\vec\ga\land_\fI\vec\gb\in F$ is the supremum in $\fI$ of the cone generated by the sequence $(\max(\ga_j,\gb_j))_{j\in\gw}$. Since $\ga_{i+1}\leq \ell(\ga_i)<\gb_{i+1}$ we have $\gy_{i+1}\geq \max(\ga_{i+1},\gb_{i+1})=\gb_{i+1}$. We have $\gy_i\geq \ga_i$, so one can write $\gy_i=\ga_i+\gd$, for some $\gd$. Moreover, $\ell(\gy_i)\geq \gy_{i+1}\geq \gb_{i+1}$ whereas $\ell(\ga_i)<\gb_{i+1}$, hence $\gd>0$. It follows that $\ell(\gy_i)=\ell(\gd)$, therefore $\gd\geq \gw^{\ell(\gy_i)}\geq \gw^{\gb_{i+1}}$ and $\gy_i\geq \ga_i+\gw^{\gb_{i+1}}$. Since $\gy_i\in F_i$, by item (ii) we obtain the claim. \ep

\bl\label{f-clos-inv}
Let $(F_i)_{i\in\gw}$ be a sequence of subsets of $\ge_0$ satisfying the conditions of Lemma \ref{f-clos}. Then the set $F:=\{\vec\ga\in I:\al{i} \ga_i\in F_i\}$ is a filter in $\fI$.
 \el

\bp\ We consider two cases. If $F_i\neq \{0\}$ for infinitely many $i<\gw$, then by Condition (ii) each $F_i$ is unbounded below $\ge_0$, hence $F_i=\ge_0$, for each $i$. In this case $F$ is the improper filter.

So, we may assume that there is an $n\in\gw$ such that $F_i=\{0\}$ for all $i>n$ and $F_n\neq\{0\}$. The set $F$ is obviously upwards closed under $\leq_\fI$. We show that $\vec\ga,\vec\gb\in F$ implies $\vec\ga\land_\fI\vec\gb\in F$. Let $\vec\gy\in\cE$ be such that $\gy_i=\max(\ga_i,\gb_i)$. Notice that $\gy_i\in F_i$, for each $i\in\gw$. Recall from Lemma \ref{glb} that $\vec\gd:=\vec\ga\land_\fI\vec\gb$ is defined by the formula
$$\gd_i:= \begin{cases}
\gy_i, & \text{if $\ell(\gy_i)\geq \gd_{i+1}$,} \\
\gy_i+\gw^{\gd_{i+1}}, & \text{otherwise}.
\end{cases}
$$
Using downward induction on $i\leq n$ and Condition (ii) we see that each $\gd_i\in F_i$. \ep

With each filter $F$ in $\fI$ we associate a sequence $\vec\ga_F\in\cE$ by letting $\ga_i:=\sup\{\gb_i+1:\vec\gb\in F\}$, for each $i\in\gw$. (For obvious reasons, $\ga_i>0$.)

\bl \label{inab} For each filter $F$, the sequence $\vec\ga=\vec\ga_F$ satisfies the following condition: For all $i\in\gw$, either $\ga_i$ is a limit ordinal and $\ga_{i+1}\leq \ell(\ga_i)$, or $\ga_i=\ga_i'+1$ and $\ga_{i+1}\leq \ell(\ga_i')+1$, for some $\ga_i'$.
\el

\bp\ Suppose $\ga_i=\ga_i'+1$ and show that $\ga_{i+1}\leq\ell(\ga_i')+1$. This means that $F_i=[0,\ga_i']$ and we need to show that $F_{i+1}\subseteq [0,\ell(\ga_i')]$. Assume otherwise and pick a $\gy\in F_{i+1}$ such that $\gy>\ell(\ga_i')$. Then, by Lemma \ref{f-clos} (ii), $\ga_i'+\gw^{\gy}\in F_i$ contradicting $F_i\subseteq [0,\ga_i']$.

Now suppose $\ga_i$ is a limit ordinal. Then $F_i=[0,\ga_i)$ and we need to prove that $F_{i+1}\subseteq [0,\ell(\ga_i))$.
Assume $\gy>\ell(\ga_{i})$ and $\gy\in F_{i+1}$. The ordinal $\ga_i$ can be represented as $\ga_i=\gb+\gw^{\ell(\ga_i)}$, for some $\gb<\ga_i$. Then $\gb+1<\ga_i$ and by Lemma \ref{f-clos} (ii), $\gb+1+\gw^{\gy}\in F_i$. However, $\gw^\gy>\gw^{\ell(\ga_i)}$ and hence $\gb+\gw^{\gy}>\gb+\gw^{\ell(\ga_i)}=\ga_i$, a contradiction. \ep

A sequence $\vec\ga\in\cE$ satisfying the condition of Lemma \ref{inab} will be called \emph{suitable}. For each $\vec\ga\in\cE$ define $F_{\vec\ga}:=\{\vec\gb\in I:\al{i\in\gw} \ga_i>\gb_i\}.$

\bl\ \label{suit} Suppose an $\vec\ga\in\cE$ is suitable. Then $F_{\vec\ga}$ is a filter in $\fI$.
\el
\bp\ Using Lemma \ref{f-clos-inv}, it is sufficient to prove for each $i\in\gw$ that, if $\ga<\ga_i$, $\gb<\ga_{i+1}$ and $\ell(\ga)<\gb$ then $\ga+\gw^\gb<\ga_i$. (The other condition is clearly satisfied.) We consider two cases.

If $\ga_i=\ga_i'+1$ then we have $\ga\leq\ga_i'$ and $\gb\leq\ell(\ga_i')$, since $\gb<\ga_{i+1}\leq\ell(\ga_i')+1$. Since we assume $\ell(\ga)<\gb$, we must have $\ga\neq \ga_i'$ and hence $\ga<\ga_i'$. But then $\ga+\gw^{\gb}\leq \ga+\gw^{\ell(\ga_i')}\leq \ga_i'$.

If $\ga_i$ is a limit ordinal, then $\gb<\ga_{i+1}\leq\ell(\ga_i)$. Since $\ga<\ga_i$ one can find a $\gy>0$ such that $\ga+\gy=\ga_i$. Then $\ell(\ga_i)=\ell(\gy)$, hence $\gw^\gb<\gw^{\ell(\ga_i)}\leq \gy$. It follows that $\ga+\gw^{\gb}<\ga+\gy=\ga_i$. \ep

Let $I^\nab$ denote the set of all suitable sequences.
\bl The maps $F\mapsto\vec\ga_F$ and $\vec\ga\mapsto F_{\vec\ga}$ are mutually inverse  bijections between $\cF(\fI)$ and $I^\nab$.
\el
\bp\ To prove $F_{\vec\ga_F}=F$ we must show that $F=\{\vec\gb\in I:\al{i} \gb_i<\ga_i\}$ where $\vec\ga=\vec\ga_F$. The inclusion $(\subseteq)$ is obvious. For the opposite inclusion  consider any $\vec\gb$ such that $\al{i} \gb_i<\ga_i$ and $\vec\gb\in I$. Let $\gb_0,\dots,\gb_n$ be all the nonzero coordinates of $\vec\gb$. Since $\gb_i<\ga_i$, for each $i\leq n$ there is a $\vec\gy^i\in F$ such that $\gy^i_i\geq \gb_i$. Let $\vec\gd$ denote the g.l.b.\ in $\fI$ of $\vec\gy^i$, for all $i\leq n$. Then $\vec\gd\in F$ and $\vec\gd\leq_\fI \vec\gb$, hence $\vec\gb\in F$.

To show that $\vec\ga_{F_{\vec\ga}}=\vec\ga$ we must prove that, for all $\vec\ga\in I^\nab$,
\beq \label{eqsup} \al{i\in\gw} \sup\{\gb_i+1:\vec\gb\in I,\ \al{j} \ga_j>\gb_j\} = \ga_i.\eeq
The inequality $(\leq)$ is obvious. For the opposite inequality consider any $\gy<\ga_i$. Let $\vec\gy:=(\gw_i(\gy),\dots,\gw^\gy,\gy,0,\dots)\in I$. By downward induction on $j\leq n$ we show that $\gy_j<\ga_j$ (for $j>n$ this is obvious since $\ga_j>0$). We assume the claim holds for a $j\leq n$ such that $j>0$ and prove it for $j-1$. Notice that by Lemmas \ref{suit} and \ref{f-clos}, if $\ga<\ga_{j-1}$, $\gb<\ga_j$ and $\ell(\ga)<\gb$ then $\ga+\gw^\gb<\ga_{j-1}$. Take $0$ for $\ga$ and $\gy_j$ for $\gb$ and conclude $\gy_{j-1}=\gw^{\gy_j}<\ga_{j-1}$ which proves the induction step.

Since $\gy_i=\gy$, $\vec\gy\in I$ and $\al{j\in\gw}\gy_j<\ga_j$, we obtain $\gy<\sup\{\gb_i+1:\vec\gb\in I,\ \al{j}\: \ga_j>\gb_j\}$. Since $\gy<\ga_i$ was arbitrary, this yields the claim \refeq{eqsup}. \ep

Thus, we have effectively characterized the domain of the universal Kripke frame $\cF(\fI)$. The next task is to characterize the relations $R_n$ and $S_n$ on $\cF(\fI)$ corresponding to $\Diamond_n$ and $\nab_n$, for all $n\in\gw$, respectively. We have, for all filters $F,G\in \cF(\fI)$,
\begin{eqnarray*}
F R_n G & \iff & \al{\vec\gy\in G} \Diamond^\fI_n\vec\gy \in F; \\
F S_n G & \iff & \al{\vec\gy\in G} \nab_n^\fI\vec\gy \in F.
\end{eqnarray*}

Let $\vec\ga:=\vec\ga_G$ and $\vec\gb:=\vec\ga_F$. Recall that $\nab_n\vec\gy$ is the truncation of $\vec\gy$ at position $n$ (followed by zeros). Let $G\rst n$ denote the set of all such truncated sequences, for all $\vec\gy\in G$.  We obviously have:
\begin{eqnarray*}
F S_n G & \iff & G\rst n \subseteq F \\
 & \iff & \al{i\leq n} \ga_i\leq \gb_i.
\end{eqnarray*}

A characterization of $R_n$ is more complicated. We begin with an auxiliary notion. Given a filter $G$ in $\fI$ let $\Diamond_n^\fI G$ denote the upwards closed subset of $\fI$ generated by $\{\Diamond_n^\fI\vec\ga:\vec\ga\in G\}$.
Since $G$ is a filter and $\Diamond_n^\fI(\vec\ga\land\vec\gb)\leq_\fI \Diamond_n^\fI\vec\ga \land\Diamond_n^\fI\vec\gb$ holds in $\fI$, the set $\Diamond_n^\fI G$ is also a filter in $\fI$.
 Clearly, $F R_n G$ holds iff $\Diamond_n^\fI G\subseteq F$.

\bl\ \label{ha} Let $H=\Diamond_n^\fI G$ and $G_i,H_i$ denote the $i$-th projections of $G$ and $H$, respectively. Then $H_{n+1}=\{0\}$ and, for all $i\leq n$, $H_i$ is the initial segment generated by
$$\{\gy+\gw^\gd:\gy\in G_i, \gd\in H_{i+1}\}.$$
\el
\bp\ The claim $H_{n+1}=\{0\}$ is obvious. We denote the set $\{\gy+\gw^\gd:\gy\in G_i, \gd\in H_{i+1}\}$ by $H'_i$ and show that $H_i$ is the initial segment generated by $H'_i$, for all $i\leq n$.

First, we prove $H'_i\subseteq H_i$. Assume $i\leq n$, $\gy\in G_i$ and $\gd\in H_{i+1}$. Then there are $\vec\ga\in G$ and $\vec\gb\in H$ such that $\ga_i=\gy$ and $\gb_{i+1}=\gd$. We have that $\vec\gb\geq_\fI\Diamond_n^\fI \vec\gb'$, for some $\vec\gb'\in G$. Let $\vec\ga':=\vec\ga\land_\fI \vec\gb'\in G$ and let $\vec\ga'':=\Diamond_n^\fI \vec\ga'\in H$.

Since $\vec\ga'\leq_\fI \vec\gb'$ we have $\Diamond_n^\fI \vec\ga'\leq_\fI \Diamond_n^\fI\vec\gb'$, hence $\ga''_{i+1}\geq  \gb_{i+1}=\gd$. Also, $\ga_i'\geq\ga_i=\gy$. It follows that $\ga''_i=\ga'_i+\gw^{\ga_{i+1}''}\geq \gy+\gw^\gd$. Therefore, $\gy+\gw^\gd\in H_i$.

Second, we prove $H_i$ is contained in the initial segment generated by $H'_i$. Suppose $\ga\in H_i$, then there is a $\vec\gb\in G$ such that $\gb'_i\geq \ga$ where $\vec\gb':=\Diamond_n^\fI\vec\gb$. Since $\vec\gb'\in \Diamond_n^\fI G = H$ we clearly have $\gb'_{i+1}\in H_{i+1}$.
By the definition of $\Diamond_n^\fI$, $\gb'_i:= \gb_i+\gw^{\gb'_{i+1}}$ where $\gb_i\in G_i$ and $\gb'_{i+1}\in H_{i+1}$. Hence, $\gb'_i\in H'_i$ and $\ga\leq \gb'_i$ is contained in the initial segment generated by $H'_i$.
\ep

Let $\vec\ga:=\vec\ga_G$ and let $\vec\nu:=\vec\ga_H$ denote the sequence corresponding to the filter $\Diamond_n^\fI G$. From Lemma~\ref{ha} we conclude that
$\nu_i=1$, for all $i>n$, and otherwise $$\nu_i:=\sup\{\gy+\gw^\gd+1:\gy<\ga_{i},\gd<\nu_{i+1}\}.$$
We are going to give an explicit formula for calculating $\nu_i$.
Given a sequence $\vec\ga\in I^\nab$ define a sequence $\gs_n(\vec\ga):=(\ga'_0,\ga'_1,\ga'_2,\dots,\ga'_n,1,1,\dots)$ where $\ga'_{n+1}=1$
\ignore{
\[
\ga'_n:=\begin{cases} \ga_n, \quad \text{if $\ga_n\in\Lim$}; \\
\ga_n+1, \quad \text{otherwise}.
\end{cases}
\]
}
and recursively, for all $i=n,n-1,\dots,0$,
\[
\ga'_i:=\begin{cases} \ga_i, \quad \text{if $\ell(\ga_i)\geq \ga'_{i+1}$; otherwise} \\ \ga_i+\gw^{\ga'_{i+1}}, \quad \text{if $\ga'_{i+1}\in\Lim$;} \\
 \ga_i+\gw^{\gd}+1, \quad \text{if $\ga'_{i+1}=\gd+1$ and $\ell(\ga_i)<\gd$;} \\
 \ga_i+1, \quad \text{if $\ga'_{i+1}=\gd+1$ and $\ell(\ga_i)=\gd$}.
\end{cases}
\]

\bl\ For each $\vec\ga\in I^\nab$, $\vec\nu=\gs_n(\vec\ga)$.
\el
\bp\ We are going to show that, for each $i\leq n+1$, $\nu_i=\ga'_i$.

For $i=n$, $\ga'_n=\ga_n$ if $\ga_n\in \Lim$, otherwise $\ga'_n=\ga_n+1$. In each case,
$\ga'_n=\sup\{\gy+1+1:\gy<\ga_n\}=\nu_n$ and the claim holds.

Consider the case $i<n$. By the induction hypothesis $\ga'_{i+1}=\nu_{i+1}$. We show $$\ga'_i=\sup\{\gy+\gw^\gd+1:\gy<\ga_{i},\gd<\ga'_{i+1}\}=\nu_i,$$
where the second equality holds by the induction hypothesis.
\begin{description}
\item[Case 1:] $\ell(\ga_i)\geq \ga'_{i+1}$. Then, for all $\gd<\ga'_{i+1}$ and $\gy_i<\ga_i$, $$\gy_i+\gw^\gd<\gy_i+\gw^{\ga'_{i+1}}\leq \ga_i.$$
It follows that $\nu_i\leq\ga_i=\ga'_i$. The opposite inequality $\ga_i\leq\nu_{i}$ is obvious.

\item[Case 2:] $\ell(\ga_i)< \ga'_{i+1}$. We consider two subcases.

\begin{description} \item[Subcase 2.1:] $\ga'_{i+1}\in \Lim$. Then $\nu_i=\ga_i+\gw^{\ga'_{i+1}}=\ga'_i$.

\item[Subcase 2.2:] $\ga'_{i+1}\notin\Lim$. Then, for some $\gd$, $\ga'_{i+1}=\gd+1$.

 \bi \item If $\ell(\ga_i)=\gd$ then $\ga_i=\gy+\gw^\gd$, for some $\gy<\ga_i$. Hence, $\nu_i\geq \gy+\gw^\gd+1=\ga_i+1=\ga'_i$.

 On the other hand, if $\gy_i<\ga_i$ then either $\gy_i<\gy$ or $\gy_i=\gy+\ge$ with $\ge<\gw^\gd$. In each case $\gy_i+\gw^\gd\leq \gy+\gw^\gd=\ga_i$. Hence, $\gy_i+\gw^\gd+1\leq\ga'_i$ and $\nu_i\leq\ga'_i$.

\item If $\ell(\ga_i)<\gd$ then $\ga_i=\gy+\gw^\ge$, for some $\gy<\ga_i$ and some $\ge<\gd$. Hence, $\nu_i\geq \gy+\gw^\gd+1=\gy+\gw^\ge+\gw^\gd+1=\ga'_i$. On the other hand, if $\gy_i<\ga_i$ then $\gy_i+\gw^\gd+1\leq \ga_i+\gw^\gd+1=\ga'_i$. Hence $\nu_i=\ga'_i$.
    \ei
\end{description}
\end{description}
\ep
Now we obtain the desired characterization of relations $R_n$ on $\fI^\nab$. \\ Let $\vec\ga:= \vec\ga_G$ and $\vec\gb:=\vec\ga_F$.
\bcor For each $n\in\gw$,
$F R_n G$ iff $\vec\gb\leq_\fI \gs_n(\vec\ga)$ iff $\al{i\leq n} \ga'_i\leq \gb_i$.
\ecor

\bibliographystyle{plain}
\bibliography{ref-all2}
\end{document}